\documentclass[11pt, reqno]{amsart}

\usepackage{amsmath}
\usepackage{fullpage}
\usepackage{amssymb}
\usepackage{amsthm}
\usepackage{bbm}
\usepackage{subfigure}
\usepackage{float}
\usepackage{graphicx}
\usepackage{color}

\newtheorem{lemma}{Lemma}[section]
\newtheorem{theorem}[lemma]{Theorem}
\newtheorem{proposition}[lemma]{Proposition}
\newtheorem{remark}[lemma]{Remark}

\graphicspath{{./figs/}}

\title{Low regularity error estimates for\\ high dimensional nonlinear Schr\"{o}dinger equations}

\author{Lun Ji}
\address{Department of Mathematics, Universit\"{a}t Innsbruck, Technikerstr.~13, 6020 Innsbruck, Austria (L.~Ji)}
\email{lun.ji@uibk.ac.at}

\author{Alexander Ostermann}
\address{Department of Mathematics, Universit\"{a}t Innsbruck, Technikerstr.~13, 6020 Innsbruck, Austria (A.~Ostermann)}
\email{alexander.ostermann@uibk.ac.at}

\begin{document}

\begin{abstract}
The filtered Lie splitting scheme is an established method for the numerical integration of the periodic nonlinear Schr\"{o}dinger equation at low regularity. Its temporal convergence was recently analyzed in a framework of discrete Bourgain spaces in one and two space dimensions for initial data in $H^s$ with $0<s\leq 2$. Here, this analysis is extended to dimensions $d=3, 4, 5$ for data satisfying $d/2-1 < s \leq 2$. In this setting, convergence of order $s/2$ in $L^2$ is proven. Numerical examples illustrate these convergence results.
\end{abstract}

\maketitle

\section{Introduction}

The convergence analysis of numerical discretizations of the periodic nonlinear Schr\"{o}dinger equation is typically carried out in a Sobolev space setting of periodic functions, where the nonlinearity is treated with standard Lipschitz techniques~\cite{Bao, Besse, Lubich08}. For example, the cubic nonlinear Schr\"{o}dinger equation (NLS)
\begin{equation}\label{nls}
i\partial_t u=-\Delta u-\mu|u|^2u,\quad (t,x)\in\mathbb{R}\times\mathbb{T}^d,\quad \mu\in\{\pm 1\}
\end{equation}
on the $d$-dimensional torus $\mathbb{T}^d$ with initial data $u(0)\in H^s(\mathbb{T}^d)$ for $s > d/2$ can be handled in this way. In order to minimize the additional spatial regularity assumptions that arise for standard integration schemes~\cite{Lubich08}, problem adapted integrators such as Fourier integrators have been constructed in the literature. For details, we refer to \cite{ESS16, KOS19, OsWY, OsY, OS18} and references therein.

The case of less regular initial data with $s\leq d/2$, however, is much more involved. The usual framework of Sobolev spaces is not applicable due to the lack of suitable embedding results. In addition, low and high frequencies have to be treated differently, and special filtering techniques for the integrators become essential~\cite{Ignat}. All this requires a new framework for the analysis. Building on the analytical results of Bourgain \cite{Bourg}, we have constructed discrete Bourgain spaces and proved the convergence of filtered integrators for \eqref{nls} in space dimensions one and two, see~\cite{JiORS,Ost,Ost1}. In particular, we considered in these papers the filtered Lie splitting scheme
\begin{equation}\label{lie}
u_{n+1}=\Psi^\tau(u_n)=e^{i\tau\Delta}\Pi_\tau(e^{\mu i\tau|\Pi_\tau u_n|^2}\Pi_\tau u_n),\quad u_0=\Pi_\tau u(0),
\end{equation}
where the projection operator $\Pi_\tau$ for $\tau>0$ is defined as:
\begin{equation}\label{proj}
\Pi_\tau=\overline{\Pi}_\tau=\chi\left(\dfrac{-i\nabla}{\tau^{-\frac{1}{2}}}\right).
\end{equation}
Here $\chi$ denotes the characteristic function of the unit cube $[-1,1]^d$. Similar filtering techniques are widely used in the analysis of low regularity problems \cite{Ignat, Rou, Wu}. Note that our filter depends on the size of the time step $\tau$ and introduces a CFL condition into the scheme by the back door.

The convergence of filtered Lie splitting for $0<s\leq2$ was first analyzed in dimension one in \cite{Ost1}, where convergence of order $\frac{s}2$ was proved. This result was later improved by using more involved filtering techniques in \cite{Wu}, and extended to dimension two in \cite{JiORS}.

Throughout the paper, we denote by $s_0$, $b_0$ and $b_1$ some real numbers satisfying the restrictions
\begin{equation}\label{s0b0}
s_0>\max(0,\tfrac d2-1) ,\quad b_0\in\big(\tfrac12,\min(\tfrac12+\tfrac14(s_0-\tfrac d2+1),\tfrac34)\big),\quad b_1=1-b_0,
\end{equation}
where $d\geq1$ denotes the space dimension. 

The purpose of this paper is to give a general convergence proof for method (2) in dimensions 1 to~5. Many ideas for such a proof can be taken from our previous proofs for the cases $d = 1$ and $d = 2$, but several estimates indeed depend on the dimension. In order to give a readable proof in general dimension, we have chosen to illustrate the local error analysis in dimension~3. However, unlike in our works~\cite{JiORS,Ost1},  we have used dimension dependent parameters this time. They have to be chosen accordingly to obtain the corresponding proof in the desired dimension. In Table 1 we give these parameters explicitly for the range $1\le d \le 5$. The global error, however, is directly estimated in dimension~$d$.

The main achievement of our work is the following convergence result.

\begin{theorem}\label{mainthm}
For given $1\leq d\leq5$ and $T>0$, let $u_0\in H^{s_0}(\mathbb{T}^d)$ be such that the exact solution $u$ of \eqref{nls} with initial data $u_0$ exists up to time $T$ for $s_0$ as in \eqref{s0b0} satisfying $s_0\le 2$. Let $u_n$ be the numerical solution defined by the scheme \eqref{lie}. Then, we have the following error estimate: there exist $\tau_0>0$ and $C_T>0$ such that for every time step $\tau\in(0,\tau_0]$
$$
\Vert u_n-u(t_n)\Vert_{L^2(\mathbb{T}^d)}\leq C_T\tau^{\frac{s_0}{2}},\quad 0\leq n\tau\leq T.
$$
\end{theorem}

To simplify notation, we will henceforth consider the defocusing case ($\mu=-1$) only. We stress, however, that our analysis given below remains true \emph{mutatis mutandis} for $\mu=1$ as long as the solution exists up to time $T$.

\subsection*{Outline of the paper.} The paper is organized as follows. In section \ref{sectionbourgain}, we give a brief outline of continuous and discrete Bourgain spaces and we introduce the employed analytic framework. The local error of the filtered Lie splitting scheme (\ref{lie}) is studied in section \ref{sectionlocal}. Global error estimates are given in section \ref{sectionglobal}, where our convergence result, Theorem \ref{mainthm}, is proved. Numerical examples in dimensions 3 and 4 illustrate this analytical result.

\subsection*{Notations.} We write $A\lesssim B$ to express that $A\leq CB$ for some generic constant $C$. In particular, this constant must not depend on the time step size $\tau \in (0, 1]$. If we want to emphasize that the constant $C$ depends on a specific parameter $\gamma$, we use the symbol $\lesssim_\gamma$. Moreover, $A\sim B$ stands for $A\lesssim B\lesssim A$.

We also use the Japanese bracket notation $\langle\,\cdot\,\rangle=(1+|\cdot|^2)^\frac12$ and, for sequences $(u_{n})_{n \in \mathbb{Z}} \in  X^\mathbb{Z}$ in a Banach space $X$ with norm $\|\cdot \|_{X}$, we employ the standard norms
\begin{equation}\label{lpnorm}
\|u_n\|_{l^p_\tau X}=\left( \tau \sum_n \|u_n\|_X^p\right)^\frac1p, \quad \|u_n\|_{l^\infty_\tau X} = \sup_{n\in\mathbb{Z}} \|u_n\|_X.
\end{equation}

\section{A Bourgain space framework}\label{sectionbourgain}

In this section, we will recall the construction of Bourgain spaces \cite{Bourg}, discrete Bourgain spaces \cite{JiORS, Ost} and some useful properties. As a first step we would recall the definition of Bourgain spaces.

For a function $u(t,x)$ defined on $\mathbb{R}\times\mathbb{T}^d$, we denote the time-space Fourier transform by $\tilde{u}(\sigma,k)$, i.e.
$$
\tilde{u}(\sigma,k)=\int_{\mathbb{R}\times\mathbb{T}^d}u(t,x)e^{-i\sigma t-i\langle k,x\rangle}dx dt,\qquad \sigma\in\mathbb{R}, k\in\mathbb{Z}^d,
$$
where $\langle\cdot,\cdot\rangle$ is the inner product in $\mathbb R^d$. The inverse transform is given by
$$
u(t,x)=\frac1{(2\pi)^d}\sum\limits_{k\in\mathbb{Z}^d}\hat{u}_k(t)e^{i\langle k,x\rangle},
$$
where the Fourier coefficients are given as $\hat{u}_k(t)=\frac1{2\pi}\int_{\mathbb{R}}\tilde{u}(\sigma,k)e^{i\sigma t}d\sigma$.

The Bourgain space $X^{s,b}= X^{s,b}(\mathbb{R}\times\mathbb{T}^d)$, as a Banach space, is defined by the norm
$$
\Vert u\Vert_{X^{s,b}}=\Vert\langle k\rangle^s\langle\sigma + |k|^2\rangle^b\tilde{u}(\sigma,k)\Vert_{L^2 l^2}
$$
and contains all functions with finite norm. Note that, in particular, $X^{s,0}=L^2H^s$.

Then, we recall from \cite{Tao} some useful properties of Bourgain spaces.
\begin{lemma}\label{contiprop}
For any $\eta\in\mathcal{C}_c^\infty(\mathbb{R})$, $s\in\mathbb{R}$ and functions in the corresponding spaces, we have that
\begin{align}
\label{iniest}\Vert\eta(t)e^{it\Delta}f\Vert_{X^{s,b}}&\lesssim_{\eta,b}\Vert f\Vert_{H^s(\mathbb{T}^d)},\quad b\in\mathbb{R},\\
\label{d}\Vert u\Vert_{L^\infty H^s}&\lesssim_b\Vert u\Vert_{X^{s,b}},\quad b>\tfrac12\\
\label{b-1}\left\Vert\eta(t)\int_{-\infty}^te^{i(t-t')\Delta} F(t')dt'\right\Vert_{X^{s,b}(T)}&\lesssim_{\eta,b}\Vert F\Vert_{X^{s,b-1}(T)},\quad b>\tfrac12.
\end{align}
\end{lemma}

These estimates are dimension independent and proven, for example, in \cite[section 2.6]{Tao}.

We point out that method \eqref{lie} can also be seen as standard Lie splitting, applied to the filtered equation
\begin{equation}\label{projeq}
i\partial_t u^\tau=-\Delta u^\tau-\mu \Pi_\tau(|\Pi_\tau u^\tau|^2\Pi_\tau u^\tau), \quad u^\tau(0)=\Pi_\tau u(0).
\end{equation}
We will make heavily use of this property.

The next step is to give some important properties related to the solutions $u$ and $u^\tau$ to the equations \eqref{nls} and \eqref{projeq}, respectively.

\begin{theorem}\label{theoutau}
For $s_0,~b_0$ defined in \eqref{s0b0}, let $u_0\in H^{s_0}$. Then, there exists $T>0$ and solutions $u$ to \eqref{nls} and $u^\tau$ to \eqref{projeq} on $[0,T]$ with initial data $u(0)=u_0$, where $u,u^\tau\in X^{s_0,b_0}$ and $\tau$ is sufficiently small. Moreover, we have that
\begin{align}
\label{boundedutau}\|u^\tau\|_{X^{s_0,b_0}}&\leq C_T,\\
\label{diffutau}\|u-u^\tau\|_{X^{0, b_0}} &\leq C_T\tau^\frac{s_0}2,\\
\label{utaumieux}\|u^\tau\|_{X^{s_0,1}}&\leq C_{T,\sigma}\tau^{-\sigma}
\end{align}
for some $C_T,C_{T,\sigma}>0$ and any $\sigma>\frac12$.
\end{theorem}

\begin{proof}
The proof of \eqref{boundedutau} and \eqref{diffutau} can be extended in a natural way from the two-dimensional case given in \cite[section 2]{JiORS}. We omit the details here. Note that both $u$ and $u^\tau$ are defined globally, but are solutions to the equations only on $[0,T]$ (see also \cite[Remark 2.7]{JiORS} and \cite[Proposition 4.3]{Ost1}).

For \eqref{utaumieux}, by Duhamel's formula, we have that $u^\tau$ solves the equation (see \cite{JiORS,Ost1})
$$
u^{\tau}(t)=e^{it\Delta}u_0+\int_0^te^{i(t-\vartheta)\Delta} F(u^{\tau}(\vartheta))d\vartheta, \quad F(v)=-i\Pi_\tau(|\Pi_\tau v|^2\Pi_\tau v)
$$
on $[0,T]$. Using \eqref{iniest} and \eqref{b-1}, we first get
$$
\|u^\tau\|_{X^{s_0,1}}\lesssim\|u_{0}\|_{H^{s_0}}+\Vert F(u^{\tau})\Vert_{L^2H^{s_0}}=\|u_0\|_{H^{s_0}}+\Vert\langle\partial_x\rangle^{s_0}F(u^\tau)\Vert_{L^2}.
$$
By the Kato--Ponce inequality (see \cite{BourgL, Kato, Mus}), we obtain
$$
\Vert\langle\partial_x\rangle^{s_0}F(u^\tau)\Vert_{L^2}\lesssim \Vert\langle\partial_x\rangle^{s_0}\Pi_\tau u^{\tau}\Vert_{L^p}\Vert\Pi_\tau u^{\tau}\Vert^2_{L^q},
$$
where $p=2+\frac4d$ and $q=2d+4$. Note that $\frac1p+\frac2q=\frac12$. Thus by employing the Strichartz estimates given in \cite{Bourg}, we have for any $\varepsilon>0$ sufficiently small that
\begin{align*}
\Vert\langle\partial_x\rangle^{s_0}\Pi_\tau u^{\tau}\Vert_{L^p}&\lesssim\Vert\langle\partial_x\rangle^{s_0}\Pi_\tau u^{\tau}\Vert_{X^{\varepsilon,b_0}}=\Vert\Pi_\tau u^{\tau}\Vert_{X^{s_0+\varepsilon,b_0}},\\
\Vert\Pi_\tau u^{\tau}\Vert_{L^q}&\lesssim\Vert\Pi_\tau u^{\tau}\Vert_{X^{\frac d2-\frac12+\varepsilon,b_0}}\leq\Vert\Pi_\tau u^{\tau}\Vert_{X^{s_0+\frac12,b_0}}
\end{align*}
since $s_0>\frac d2-1$, where $b_0$ is defined in \eqref{s0b0}. Note that $\Pi_\tau$ is projecting on spatial frequencies smaller than $\tau^{-\frac12}$. We thus get
$$
\|u^\tau\|_{X^{s_0,1}}\lesssim\|u_0\|_{H^{s_0}}+\|\Pi_\tau u^\tau\|_{X^{s_0+\varepsilon,b_0}}\|\Pi_\tau u^\tau\|^2_{X^{s_0+\frac12,b_0}}\lesssim\|u_0\|_{H^{s_0}} + \tau^{-\frac12-\frac12\varepsilon}\|u^\tau\|^3_{X^{s_0, b_0}}.
$$
Finally, by using \eqref{boundedutau}, we can conclude \eqref{utaumieux} by taking $\sigma=\frac12+\frac12\varepsilon>\frac12$, which can be chosen arbitrarily close to $\frac12$.
\end{proof}

Next, we shall recall the definition of discrete Bourgain spaces, and discuss some of their properties. For more details on the properties and the proofs, see also \cite{JiORS,Ost}.

First, we take $(u_n(x))_n$ to be a sequence of functions on the torus $\mathbb{T}^d$, with its time-space Fourier transform
$$
\widetilde{u_n}(\sigma,k)=\tau\sum\limits_{m\in\mathbb{Z}}\widehat{u_m}(k)e^{im\tau\sigma},
$$
where
$$
\widehat{u_m}(k)=\dfrac{1}{(2\pi)^d}\int_{\mathbb{T}^d}u_m(x)e^{-i\langle k,x\rangle}dx.
$$
Note that, in this framework, $\widetilde{u_{n}}$ is a $2\pi/\tau$ periodic function in $\sigma$.

Then the discrete Bourgain space $X^{s,b}_\tau$ will be defined as the set of all sequences of functions $(u_n(x))_n$ with the norm
\begin{equation}\label{5}
\Vert u_n\Vert_{X^{s,b}_\tau}=\Vert\langle k\rangle^s\langle d_\tau(\sigma-|k|^2)\rangle^b\widetilde{u_n}(\sigma,k)\Vert_{L^2l^2},
\end{equation}
where $d_\tau(\sigma)=\frac{e^{i\tau\sigma}-1}{\tau}$. For any fixed $(u_n)_n$, the norm is obviously an increasing function for both $s$ and $b$.

From definitions \eqref{proj} and \eqref{5}, we can directly get the following properties: for any $s\geqslant s^\prime$ and $b\geqslant b^\prime$, we have
\begin{align}
\label{b}\Vert\Pi_\tau u_n\Vert_{X^{s,b}_\tau}&\lesssim\tau^{b^\prime-b}\Vert\Pi_\tau u_n\Vert_{X^{s,b^\prime}_\tau},\\
\label{t}\Vert\Pi_\tau u_n\Vert_{X^{s,b}_\tau}&\lesssim\tau^{\frac{s^\prime-s}{2}}\Vert\Pi_\tau u_n\Vert_{X^{s^\prime,b}_\tau}.
\end{align}

In the next lemma, we will give an estimate, which will be useful for the local and global error analysis.
\begin{lemma}
For any $s\in\mathbb{R},~b>\tfrac12,~\tau\in(0,1]$ and $u_n\in l_\tau^1 H^s$, we have that
\begin{equation}\label{ydual}
\Vert u_n\Vert_{X^{s,-b}_\tau}\lesssim_b\Vert u_n\Vert_{l_\tau^1 H^s}.
\end{equation}
\end{lemma}

\begin{proof}
It suffices to prove $u_n\in X^{s,-b}_\tau$.

For any $\gamma\in\mathbb{R}$ and $v_n\in X^{\gamma,b}_\tau$, we have the estimate (see \cite{Ost})
\begin{equation}\label{y}
\Vert v_n\Vert_{l_\tau^\infty H^\gamma}\lesssim_b\Vert v_n\Vert_{X^{\gamma,b}_\tau},
\end{equation}
thus we get
$$
\tau\sum\limits_n\int_{\mathbb{T}^d}|u_nv_n|\,dx\leq\|u_n\|_{l_\tau^1 H^s}\Vert v_n\Vert_{l_\tau^\infty H^{-s}}\lesssim_b\|u_n\|_{l_\tau^1 H^s}\Vert v_n\Vert_{X^{-s,b}_\tau}.
$$
Since $v_n$ can be chosen arbitrarily in $X^{-s,b}_\tau$, we have $u_n$ is in the dual of $X^{-s,b}_\tau$, i.e. $u_n\in X^{s,-b}_\tau$. This ends the proof.
\end{proof}

We will conclude this section by providing key estimates for the analysis of scheme \eqref{lie}.
\begin{theorem}\label{disckeythm}
For any $s_0,~b_1$ defined in \eqref{s0b0}, we have
\begin{align}
\label{z}\Vert\Pi_\tau u_n\Vert_{l^4_\tau L^4}&\lesssim\Vert u_n\Vert_{X_\tau^{\frac {s_0}2,b_1}},\\
\label{r}\Vert\Pi_\tau(\Pi_\tau u_{n,1} \Pi_\tau \overline{u}_{n,2}\Pi_\tau u_{n,3})\Vert_{X_\tau^{s_0,-b_1}}&\lesssim\Vert u_{n,1}\Vert_{X_\tau^{s_0,b_1}}\Vert u_{n,2}\Vert_{X_\tau^{s_0,b_1}}\Vert u_{n,3}\Vert_{X_\tau^{s_0,b_1}}, \\
\label{rbis}\Vert\Pi_\tau(\Pi_\tau u_{n,1} \Pi_\tau \overline{u}_{n,2}\Pi_\tau u_{n,3})\Vert_{X_\tau^{0,-b_1}}&\lesssim\Vert u_{n, \sigma(1)}\Vert_{X_\tau^{s_0,b_1}}\Vert u_{n,\sigma(2)}\Vert_{X_\tau^{s_0,b_1}}\Vert u_{n, \sigma(3)}\Vert_{X_\tau^{0,b_1}},
\end{align}
where $(u_n)_n$, $(u_{n,i})_n$, are functions in the corresponding spaces and $\sigma$ is any permutation of $\{1, 2, 3\}$.
\end{theorem}

\begin{proof}
We first use Littlewood--Paley decompositions to prove \eqref{z}. Set
$$
\Pi_\tau u_n=\sum\limits_{A,B}u_n^{AB},
$$
where $A,B$ are dyadic numbers and $\tilde{u}_n^{AB}(\sigma,k)=\mathbbm{1}_{S_{AB}}\tilde{u}_n$, where $S_{AB}=\{(\sigma,k)\;|\; A\leq\langle k\rangle\leq2A,~B\leq\langle\sigma-|k|^2\rangle\leq2B\}$.
We shall first give two estimates of $\Vert\Pi_\tau u_n\Vert_{l_\tau^4L^4}$, which behave better with respect to $A$ and $B$ respectively. To be specific, we will prove
\begin{equation}\label{Aest}
\Vert\Pi_\tau u_n^{AB}\Vert_{l_\tau^4L^4}\lesssim A^{\frac d4}B^\frac14\Vert u_n^{AB}\Vert_{l_\tau^2L^2}
\end{equation}
and
\begin{equation}\label{Best}
\Vert\Pi_\tau u_n^{AB}\Vert_{l_\tau^4L^4}\lesssim A^{\frac{d-2}4+\frac\varepsilon2}B^\frac12\Vert u_n^{AB}\Vert_{l_\tau^2L^2},
\end{equation}
where $\varepsilon>0$ can be taken arbitrarily small.

Now we follow closely the steps and notations in \cite[Lemma 8.1]{JiORS}. To get \eqref{Aest}, it suffices to show
$$
\sum 1\lesssim A^d,
$$
where the sum extends over all $k_{ij},~i=1,2,~j=1,\ldots,d$ satisfying
\begin{equation}\label{sumset}
k_{11}+k_{21}=k_1,\ldots, k_{1d}+k_{2d}=k_d,\quad\sigma-\sum\limits_{j=1}^d(k_{1j}^2+k_{2j}^2)\in E_B.
\end{equation}
We fix $\sigma$ and $k_j$ for $1\leq j\leq d$. Since $k_{1j}+k_{2j}=k_j$ fixed, there are $d$ free variables where each of them has $\mathcal{O}(A)$ different choices. This ends the proof of \eqref{Aest}.

To get \eqref{Best}, we again follow closely the steps and notations in \cite[Lemma 8.1]{JiORS}. It thus will suffice to show
$$
\sum 1\lesssim A^{d-2+2\varepsilon}B,
$$
where again the sum extends over all $k_{ij},~i=1,2,~j=1,\ldots,d$ satisfying \eqref{sumset}.

We first choose $k_{1j},~3\leq j\leq d$ freely and thus $k_{2j},~3\leq j\leq d$ results automatically. Then the problem is a two-dimensional problem. Therefore, \eqref{Best} follows from the proof in \cite[Lemma 8.1]{JiORS}.

Next, we interpolate \eqref{Aest} and \eqref{Best}. We take $\theta\in(0,1)$ such that $2-4b_1<\theta<s_0-\frac d2+1$ and $\varepsilon<s_0-\frac d2+1-\theta$, then interpolate between \eqref{Aest} with strength $\theta$ and \eqref{Best} with strength $1-\theta$ to get
$$
\Vert\Pi_\tau u_n^{AB}\Vert_{l_\tau^4L^4}\lesssim A^{\frac{s_0}2-\delta_1}B^{b_1-\delta_2}\Vert u_n^{AB}\Vert_{l_\tau^2L^2},
$$
for some $\delta_1,\delta_2>0$. To conclude, by Cauchy--Schwarz, we have
$$
\Vert\Pi_\tau u_n\Vert_{l_\tau^4L^4}\leq\sum\limits_{A,B}\Vert\Pi_\tau u_n^{AB}\Vert_{l_\tau^4L^4}\lesssim\big(\sum\limits_{A,B}A^{-2\delta_1}B^{-2\delta_2}\big)^{\frac12}\big(\sum\limits_{A,B}A^{s_0}B^{2b_1}\Vert u_n^{AB}\Vert_{l_\tau^2L^2}^2\big)^{\frac12}\lesssim\Vert u_n\Vert_{X_\tau^{\frac{s_0}2,b_1}}.
$$

This ends the proof of \eqref{z}. By adapting the new (dimension related) parameters above to \cite[section 8]{JiORS}, one can get a proof of \eqref{r} and \eqref{rbis}. We omit the details here.
\end{proof}

\section{Local error estimates}\label{sectionlocal}

In this section, we will estimate the local error of the filtered Lie splitting method \eqref{lie}. We state the result in Theorem \ref{theolocal} and exemplify the proof for dimension $d=3$. Due to some required Sobolev embedding results, certain parameters in the proof have to be chosen dimension dependent. The actual values as a function of $d$ are given in Table~\ref{tab:1} below.

First of all, from the same computations as in \cite[section 3]{Ost1}, we can write the local error as
$$
\Psi^\tau(u^\tau(t_n))-u^\tau(t_{n+1})=ie^{i\tau\Delta}\mathcal{E}_{loc}(t_n,\tau,u^\tau),
$$
where $\Psi^\tau$ is defined in \eqref{lie} and
$$
\mathcal{E}_{loc}(t_n,\tau,u^\tau)=\mathcal{E}_1(t_n)+\mathcal{E}_2(t_n)+\mathcal{E}_3(t_n),
$$
where
\begin{align*}
\mathcal{E}_1(t_n)&=\Pi_\tau\int_0^\tau\big(e^{-i\vartheta\Delta}-1\big)\big(|\Pi_\tau u^\tau(t_n+\vartheta)|^2\Pi_\tau u^\tau(t_n+\vartheta)\big)d\vartheta,\\
\mathcal{E}_2(t_n)&=\Pi_\tau\int_0^\tau|\Pi_\tau u^\tau(t_n+\vartheta)|^2\Pi_\tau\big(u^\tau(t_n+\vartheta)-u^\tau(t_n)\big)d\vartheta,\\
\mathcal{E}_3(t_n)&=\Pi_\tau\int_0^\tau\Bigg(|\Pi_\tau u^\tau(t_n+\vartheta)|^2+\frac{e^{-i\tau|\Pi_\tau u^\tau(t_n)|^2}-1}{i\tau}\Bigg)\Pi_\tau u^\tau(t_n)d\vartheta.
\end{align*}

In the next step, we will estimate $\mathcal{E}_{loc}(t_n,\tau,u^\tau)$.

\begin{theorem}\label{theolocal}
For $1\leq d\leq 5$, $s_0\leq2$, $b_0$ and $b_1$ defined in $\eqref{s0b0}$, $u^\tau$ as in Theorem~\ref{theoutau}, we have for $\tau$ sufficiently small that
$$
\Vert\mathcal{E}_{loc}(t_n,\tau,u^\tau)\Vert_{X^{0,-b_1}_\tau}\leq C_T\tau^{1+\frac{s_0}2},
$$
where $C_T$ is uniform for $0\leq t_n\leq T$.
\end{theorem}

\begin{proof}
We illustrate here the proof for dimension $d=3$. The extension to $d=1,2,4,5$ is detailed in Remark~\ref{remlocal}.

First, we need to prove the boundedness of $u^\tau(t_n)$ in the norm of $X^{s,b}_\tau$ for suitable $s$ and $b$. In fact, combining \eqref{utaumieux} and \cite[Lemma 4.1]{JiORS}, we have
\begin{equation}\label{h}
\sup\limits_{t^\prime\in[0,4\tau]}\Vert u^\tau_n\Vert_{X_\tau^{s_0,\frac12-\varepsilon}}\leq C_T,
\end{equation}
where $\varepsilon>0$ can be arbitrarily small, and where $u^\tau_n(x)=u^\tau(n\tau+t^\prime,x)$ for $t^\prime\in[0,4\tau]$.

Then, similar to \cite[Theorem 5.1]{JiORS}, we have
\begin{equation}\label{e1}
\Vert\mathcal{E}_1(t_n)\Vert_{X^{0,-b_1}_\tau}\lesssim C_T\tau^{1+\frac{s_{0}}{2}}.
\end{equation}

Next, we split $\mathcal{E}_2(t_n)=\mathcal{E}_{2,1}(t_n)+\mathcal{E}_{2,2}(t_n)$, where
\begin{align*}
&\mathcal{E}_{2,1}(t_n)=\Pi_\tau\int_0^\tau|\Pi_\tau u^\tau(t_n+\vartheta)|^2\Pi_\tau\big((e^{i\vartheta\Delta}-1)u^\tau(t_n)\big)d\vartheta\\
&\mathcal{E}_{2,2}(t_n)=-i\Pi_\tau\int_0^\tau|\Pi_\tau u^\tau(t_n+\vartheta)|^2\Pi_\tau\int_0^\vartheta e^{i(\vartheta-\xi)\Delta}\Pi_\tau\big(|\Pi_\tau u^\tau(t_n+\xi)|^2\Pi_\tau u^\tau(t_n+\xi)\big)d\xi d\vartheta.
\end{align*}
Again by \cite[Theorem 5.1]{JiORS}, we have
\begin{equation}\label{e21}
\Vert\mathcal{E}_{2,1}(t_n)\Vert_{X^{0,-b_1}_\tau}\lesssim C_T\tau^{1+\frac{s_0}2}.
\end{equation}
By using \eqref{b} and \eqref{ydual}, we get that
\begin{equation}\label{e22st}
\Vert\mathcal{E}_{2,2}(t_n)\Vert_{X^{0,-b_1}_\tau}\lesssim\tau^{-(b_0-\frac12)-\varepsilon}\Vert\mathcal{E}_{2,2}(t_n)\Vert_{X^{0,-\frac12-\varepsilon}_\tau}\lesssim\tau^{2-(b_0-\frac12)-\varepsilon}\sup_{\vartheta\in[0,\tau]}\Vert (\Pi_\tau u^\tau(t_n+\vartheta))^5\Vert_{l_\tau^1L^2},
\end{equation}
for any $\varepsilon>0$ sufficiently small. By H\"older's inequality, we then get that
\begin{equation}\label{e22c0}
\Vert (\Pi_\tau u^\tau(t_n+\vartheta))^5\Vert_{l_\tau^1L^2}\lesssim\Vert \Pi_\tau u^\tau(t_n+\vartheta)\Vert^4_{l_\tau^4L^p}\Vert\Pi_\tau u^\tau(t_n+\vartheta)\Vert_{l_\tau^\infty L^q},
\end{equation}
where $\frac4p+\frac1q=\frac12$. We now split the interval $(\frac12,2]$ into three parts.

For the first case, we consider $s_0\in(\frac12,\frac45]$. By taking $p=15$ and $q=\frac{30}{7}$ in \eqref{e22c0}, we get by Sobolev embeddings, \eqref{z}, and \eqref{t} that
\begin{equation}\label{e22l4}
\begin{aligned}
\Vert\Pi_\tau u^\tau(t_n+\vartheta)\Vert_{l_\tau^4L^{15}}\lesssim\Vert\Pi_\tau u^\tau(t_n+\vartheta)\Vert_{l_\tau^4W^{\frac{11}{20},4}}&\lesssim\Vert\langle\partial_x \rangle^\frac{11}{20}\Pi_\tau u^\tau(t_n+\vartheta)\Vert_{X_\tau^{\frac14+5\varepsilon,\frac12-\varepsilon}}\\
&\lesssim\tau^{-\frac25-\frac{5\varepsilon}2+\frac{s_0}{2}}\Vert u^\tau(t_n+\vartheta)\Vert_{X_\tau^{s_0,\frac12-\varepsilon}}.
\end{aligned}
\end{equation}
Moreover, by Sobolev embeddings, the continuity of $u^\tau$ and \eqref{d}, we get that
\begin{equation}\label{e22li}
\Vert\Pi_\tau u^\tau(t_n+\vartheta)\Vert_{l_\tau^\infty L^\frac{30}7}\lesssim\Vert\Pi_\tau u^\tau(t_n+\vartheta)\Vert_{l_\tau^\infty H^\frac45}\lesssim\tau^{-\frac25+\frac{s_0}2}\Vert u^\tau\Vert_{L^\infty H^{s_0}}\lesssim\tau^{-\frac25+\frac{s_0}2}\Vert u^\tau\Vert_{X^{s_0,\frac12+\varepsilon}}.
\end{equation}
Combining \eqref{e22st}, \eqref{e22c0}, \eqref{e22l4} and \eqref{e22li}, taking $\varepsilon$ small enough ($s_0-\frac12>100\varepsilon$, for example) and using \eqref{h} and \eqref{boundedutau}, we have
\begin{equation}\label{e22c1}
\Vert\mathcal{E}_{2,2}(t_n)\Vert_{X^{0,-b_1}_\tau}\lesssim C_T\tau^{\frac{5s_0}{2}-(b_0-\frac12)-11\varepsilon}\lesssim C_T\tau^{1+\frac{s_0}{2}}
\end{equation}
since $s_0>\frac12$ and $b_0-\frac12<\frac14(s_0-\frac12)$.

For the second case, $s_0\in(\frac45,\frac32]$, we can get by a similar argument as above that
\begin{equation}\label{e22c2}
\Vert\mathcal{E}_{2,2}(t_n)\Vert_{X^{0,-b_1}_\tau}\lesssim\tau^{2-(b_0-\frac12)-\varepsilon}\Vert u^\tau(t_n+\vartheta)\Vert^4_{X_\tau^{s_0,\frac12-\varepsilon}}\Vert u^\tau\Vert_{X^{s_0,\frac12+\varepsilon}}\lesssim C_T\tau^{1+\frac{s_0}2}
\end{equation}
since $b_0-\frac12<\frac14(s_0-\frac12)$.

For the third case, $s_0\in(\frac32,2]$, by H\"older's inequality, we get the crude estimate
\begin{equation}\label{e22crude}
\Vert\mathcal{E}_{2,2}(t_n)\Vert_{X^{0,-b_1}_\tau}\lesssim\tau^2\sup_{\vartheta\in[0,\tau]}\Vert(u^\tau(t_n+\vartheta))^5\Vert_{l_\tau^2L^2}\lesssim\tau^2\sup_{\vartheta\in[0,\tau]}\Vert u^\tau(t_n+\vartheta)\Vert^2_{l_\tau^4L^p}\Vert u^\tau(t_n+\vartheta)\Vert^3_{l_\tau^\infty L^q},
\end{equation}
where $\frac2p+\frac3q=\frac12$. Taking $p=\infty$ and $q=6$, by Sobolev embeddings and again a similar argument as used in the reasoning \eqref{e22l4}--\eqref{e22c1}, we have
\begin{equation}\label{e22c3}
\Vert\mathcal{E}_{2,2}(t_n)\Vert_{X^{0,-b_1}_\tau}\lesssim\tau^2\sup_{\vartheta\in[0,\tau]}\Vert u^\tau(t_n+\vartheta)\Vert^2_{X_\tau^{s_0,\frac12-\varepsilon}}\Vert u^\tau(t_n+\vartheta)\Vert^3_{X^{s_0,\frac12+\varepsilon}}\lesssim C_T\tau^2\lesssim C_T\tau^{1+\frac{s_0}2}.
\end{equation}
Collecting \eqref{e22c1}, \eqref{e22c2} and \eqref{e22c3}, we get
\begin{equation}\label{e22}
\Vert\mathcal{E}_{2,2}(t_n)\Vert_{X^{0,-b_1}_\tau}\lesssim C_T\tau^{1+\frac{s_0}2}.
\end{equation}

To estimate $\mathcal{E}_3$, we rewrite it as
$$
\mathcal{E}_3(t_n)=\mathcal{E}_{3,1}(t_n)+\mathcal{E}_{3,2}(t_n)
$$
with
\begin{align*}
\mathcal{E}_{3,1}(t_n)&=\Pi_\tau\int_0^\tau\big(|\Pi_\tau u^\tau(t_n+\vartheta)|^2-|\Pi_\tau u^\tau(t_n)|^2\big)\Pi_\tau u^\tau(t_n)d\vartheta\\
&=\frac12\Pi_\tau\int_0^\tau\big(\Pi_\tau(\bar{u}^\tau(t_n+\vartheta)-\bar{u}^\tau(t_n))\big)\big(\Pi_\tau u^\tau(t_n+\vartheta)+\Pi_\tau u^\tau(t_n)\big)\Pi_\tau u^\tau(t_n)d\vartheta\\
&\quad+\frac12\Pi_\tau\int_0^\tau\big(\Pi_\tau(u^\tau(t_n+\vartheta)-u^\tau(t_n))\big)\big(\Pi_\tau \bar{u}^\tau(t_n+\vartheta)+\Pi_\tau \bar{u}^\tau(t_n)\big)\Pi_\tau u^\tau(t_n)d\vartheta\\
\mathcal{E}_{3,2}(t_n)&=\tau\Pi_\tau\Bigg(\dfrac{e^{-i\tau|\Pi_\tau u^\tau(t_n)|^2}-1+i\tau|\Pi_\tau u^\tau(t_n)|^2}{i\tau}\Pi_\tau u^\tau(t_n)\Bigg).
\end{align*}

Using similar arguments as in the estimate of $\mathcal{E}_2$, we can prove that
\begin{equation}\label{e31}
\Vert\mathcal{E}_{3,1}(t_n)\Vert_{X^{0,-b_1}_\tau}\lesssim C_T\tau^{1+\frac{s_{0}}{2}}.
\end{equation}

The estimate of  $\mathcal{E}_{3,2}$ is also similar since
\begin{equation}\label{expest}
\Big|\dfrac{e^{-i\tau\alpha}-1+i\tau\alpha}{i\tau}\Big|\lesssim\tau|\alpha|^2,\quad \alpha\in\mathbb{R}.
\end{equation}
Therefore, by \eqref{b} and \eqref{ydual}, we have
$$
\Vert\mathcal{E}_{3,2}(t_n)\Vert_{X^{0,-b_1}_\tau}\lesssim\tau^{-(b_0-\frac12)-\varepsilon}\Vert\mathcal{E}_{3,2}(t_n)\Vert_{X^{0,-\frac12-\varepsilon}_\tau}\lesssim\tau^{2-(b_0-\frac12)-\varepsilon}\Vert (\Pi_\tau u^\tau(t_n))^5\Vert_{l_\tau^1L^2},
$$
where $\varepsilon>0$ can be taken arbitrarily small. Consequently, by using \eqref{e22c1}, \eqref{e22c2}, we also obtain that
$$
\Vert\mathcal{E}_{3,2}(t_n)\Vert_{X^{0,-b_1}_\tau}\lesssim C_T\tau^{1+\frac{s_0}2}
$$
for the first and second cases $s_0\in(\frac12,\frac45]$ and $s_0\in(\frac45,\frac32]$.

For the third case $s_0 >\frac32$, we just use that
$$
\Vert\mathcal{E}_{3,2}(t_n)\Vert_{X^{0,-b_1}_\tau}\lesssim \Vert\mathcal{E}_{3,2}(t_n)\Vert_{X^{0,0}_\tau}\lesssim \tau^2 \|(\Pi_\tau u^\tau(t_n))^5\|_{l^2_\tau L^2}
$$
thanks to \eqref{expest}. We thus have the same upper bound as in \eqref{e22c3}. Employing the above estimates. we can get that
$$
\Vert\mathcal{E}_{3,2}(t_n)\Vert_{X^{0,-b_1}_\tau}\lesssim  C_T\tau^{1+\frac{s_0}2}
$$
for $\frac12<s_0\leq 2$. Therefore, the estimate
\begin{equation}\label{e32}
\Vert\mathcal{E}_{3,2}(t_n)\Vert_{X^{0,-b_1}_\tau}\lesssim  C_{T}\tau^{ 1 + \frac{s_0}2}
\end{equation}
holds for all $s_0\in (\frac12, 2].$

We finish the proof by collecting (\ref{e1}), (\ref{e21}), (\ref{e22}), (\ref{e31}), and (\ref{e32}).
\end{proof}

\begin{remark}\label{remlocal}
\rm{The above proof can easily be adapted to any dimension between one and five by considering appropriate intervals for $s_0$ and by carefully choosing $p$ and $q$ in \eqref{e22c0} and \eqref{e22crude}. Suitable intervals for $s_0$ and values for $p$ and $q$ are given in Table~\ref{tab:1}.

Note that we cannot expect a local error better than $\mathcal{O}(\tau^2)$, since Lie splitting is of first order. Anyhow, for $d\geq6$, using similar but slightly more involved tools, we can still get the local error bound $\mathcal{O}(\tau^2)$ for $s_0$, $b_0$ defined in \eqref{s0b0}. We omit the details here.}
\end{remark}

\begin{center}
\begin{table}[h]
\caption{Values of $s_0$, $p$ and $q$ to be chosen in the proof of Theorem~\ref{theolocal} as a function of the dimension $d$.\label{tab:1}}
\renewcommand\arraystretch{1.3333}
\begin{tabular}{||c||c|c|c|c|c||}\hline
   & $ d=1$ & $d=2$ & $d=3$ & $d=4$ & $d=5$\\ \hline\hline
$s_0$ \rm{in the first case} & $(0,\frac15] $ & $(0,\frac25]$ & $(\frac12,\frac45]$ & $(1,\frac65]$ & $(\frac32,\frac85]$ \\ \hline
$s_0$ \rm{in the second case} & $(\frac15,\frac43] $ & $(\frac25,\frac43]$ & $(\frac45,\frac32]$ & $(\frac65,\frac53]$ & $(\frac85,\frac{11}6]$ \\ \hline
$s_0$ \rm{in the third case} & $(\frac43,2]$ & $(\frac43,2]$ & $(\frac32,2]$ & $(\frac53,2]$ & $(\frac{11}6,2]$ \\ \hline
$p$ in \eqref{e22c0} & $20$ & $20$ & $15$ & $\frac{40}3$ & $\frac{25}2$ \\ \hline
$q$ in \eqref{e22c0} & $\frac{10}3$ & $\frac{10}3$ & $\frac{30}7$ & $5$ & $\frac{50}9$ \\ \hline
$p$ in \eqref{e22crude} & $\infty$ & $\infty$ & $\infty$ & $40$ & $25$ \\ \hline
$q$ in \eqref{e22crude} & $6$ & $6$ & $6$ & $\frac{20}3$ & $\frac{50}7$ \\ \hline
\end{tabular}
\end{table}
\end{center}
\section{Proof of Theorem \ref{mainthm}}\label{sectionglobal}

In this section, we will give a global error estimate and by this prove our main result, Theorem~\ref{mainthm}. We will give a proof that is valid for any dimension $2\leq d\leq5$. To get a proof for the case $d=1$, we formally set $d=2$ in the proof below.

First of all, we shall analyze the global error in the space $X_\tau^{0,b_0}$. Similar to \cite[sect.~3]{Ost1}, the global error can be written as
\begin{equation}\label{k}
e_n=-i\tau\sum\limits_{k=0}^{n-1}e^{i(n-k)\tau\Delta}(\Phi^\tau(u^\tau(t_k))-\Phi^\tau(u_k))-i\sum\limits_{k=0}^{n-1}e^{i(n-k)\tau\Delta}\mathcal{E}_{loc}(t_k,\tau,u^\tau),
\end{equation}
where
\begin{equation}\label{defphitaufin}
\Phi^\tau(w)=-\Pi_\tau\Big(\dfrac{e^{-i\tau|\Pi_\tau w|^2}-1}{i\tau}\Pi_\tau w\Big).
\end{equation}

\begin{proposition}\label{propen}
For $1\leq d\leq5,~s_0\leq2,~b_0$ and $b_1$ defined as in \eqref{s0b0}, and $\tau$ sufficiently small, we have
$$
\Vert e_n\Vert_{X^{0,b_0}_{\tau}} \lesssim \tau^\frac{s_0}2.
$$
\end{proposition}

\begin{proof}
We follow the notations in \cite[Proposition 6.1]{JiORS}. We first get from Theorem~\ref{theolocal} that
\begin{equation}\label{esten1}
\begin{aligned}
\Vert e_n\Vert_{X^{0,b_0}_\tau}\leq C_TT_1^{\varepsilon_0}\Vert \Pi_{\tau}&\bigl(F^\tau(u^\tau(t_n))-F^\tau(u^\tau(t_n)-e_n)\bigr)\Vert_{X_\tau^{0,-b'}} \\
&+ C_{T}\Vert\Pi_\tau\bigl(R^\tau(u^\tau(t_n))-R^\tau(u^\tau(t_n)-e_n) \bigr)\Vert_{X_\tau^{0,-b_1}} +C_T\tau^{\frac{s_0}2},
\end{aligned}
\end{equation}
where $\max\big(\tfrac12-\tfrac14(s_0-\tfrac d2+1),\tfrac14\big)<b^\prime<b_1$. For $F^\tau$, by \eqref{rbis} and \eqref{t}, we have the estimate
\begin{equation}\label{esten2}
\begin{aligned}
\Vert \Pi_\tau\bigl(F^\tau(u^\tau(t_n))-F^\tau(&u^\tau(t_n)-e_n) \bigr)\Vert_{X_\tau^{0,-b'}} \\
& \leq C_T \|e_n\|_{X^{0,b_1}_\tau}+C_{T,s_1}\tau^{-\frac{s_1}2}\|e_n\|_{X^{0,b_1}_\tau}^2+C_{T, s_1}\tau^{- s_1}\|e_n\|_{X^{0,b_1}_\tau}^3
\end{aligned}
\end{equation}
with $s_1$ to be chosen.

Next we estimate $\Vert\Pi_{\tau}\left(R^{\tau}(u^\tau(t_n))-R^{\tau}(u^\tau(t_n)-e_n) \right)\Vert_{X_\tau^{0,-b_1}}$. We first use \eqref{b} and \eqref{ydual} to get
\begin{align*}
\Vert \Pi_\tau\bigl(R^\tau(u^\tau(t_n))&-R^{\tau}(u^\tau(t_n)-e_n) \bigr)\Vert_{X_\tau^{0,-b_1}}\\
&\lesssim\tau^{-\delta-(\frac12-b_1)}\Vert\Pi_\tau\bigl(R^\tau (u^\tau(t_n))-R^\tau(u^\tau(t_n)-e_n)\bigr)\Vert_{X_\tau^{0,-\frac12- \delta}}\\
&\lesssim\tau^{-\delta-(\frac12-b_1)} \Vert R^{\tau}(u^\tau(t_n))-R^\tau(u^\tau(t_n)-e_n) \Vert_{l^1_\tau L^2},
\end{align*}
where $\delta>0$ to be chosen sufficiently small. Then, by the pointwise uniform in $\tau$ estimate
$$
\left| \left(R^\tau(u^\tau(t_n))-R^{\tau}(u^\tau(t_n)-e_n) \right) \right| \lesssim \tau\sum_{j=1}^5 |u^\tau(t_n)|^{5-j} |e_{n}|^j,
$$
we deduce that
\begin{equation}\label{ouf3}
\Vert\Pi_{\tau}\left(R^{\tau}(u^\tau(t_n))-R^{\tau}(u^\tau(t_n)-e_n) \right)\Vert_{X_\tau^{0,-b_1}}
\lesssim \tau^{1-\delta-(\frac12-b_1)}\sum_{j=1}^5 \left\| |\Pi_{\tau}u^\tau(t_n)|^{5-j} |\Pi_{\tau}e_n|^j \right\|_{l^1_\tau L^2}.
\end{equation}
Using H\"{o}lder's inequality, we get that
\begin{align*}
\sum_{j=1}^5 \big\| |\Pi_\tau u^\tau(&t_n)|^{5-j} |\Pi_{\tau}e_n|^j \big\|_{l^1_\tau L^2} \lesssim
\|\Pi_{\tau}u^\tau(t_n)\|_{l^4_\tau L^\infty}^4\|\Pi_{\tau} e_{n}\|_{l^\infty_\tau L^2}\\[-2.5mm]
&+ \|\Pi_{\tau}u^\tau(t_{n})\|_{l^4_\tau L^\infty}^3\|\Pi_{\tau} e_{n}\|_{l^4_\tau L^\infty}\|\Pi_{\tau} e_{n}\|_{l^\infty_\tau L^2}
+ \|\Pi_{\tau}u^\tau(t_{n})\|_{l^4_\tau L^\infty}^2\|\Pi_{\tau} e_{n}\|_{l^4_\tau L^\infty}^2\|\Pi_{\tau} e_{n}\|_{l^\infty_\tau L^2}\\
& +\|\Pi_{\tau}u^\tau(t_{n})\|_{l^4_\tau L^\infty}\|\Pi_{\tau} e_{n}\|_{l^4_\tau L^\infty}^3\|\Pi_{\tau} e_{n}\|_{l^\infty_\tau L^2} +  \|\Pi_{\tau} e_{n}\|_{l^4_\tau L^\infty}^4\|\Pi_{\tau} e_{n}\|_{l^\infty_\tau L^2}.
\end{align*}
To estimate the right-hand side, we first get by \eqref{y} that
$$
\|\Pi_\tau e_n\|_{l^\infty_\tau L^2} \lesssim \|e_n\|_{X^{0, b_0}_\tau}.
$$
Next, by Sobolev embedding and \eqref{z}, \eqref{t}, we can get that
$$
\|\Pi_{\tau} e_n\|_{l^4_\tau L^\infty} \lesssim \|e_n\|_{l^4_\tau W^{\frac d4+\delta,4}}\lesssim \|e_n\|_{X_\tau^{\frac d2-\frac12+2\delta,b_0}}\lesssim \tau^{-\frac d4+\frac14-\delta}\|e_n\|_{X_\tau^{0,b_0}}.
$$
Finally, by a similar argument, we can get that
$$
\|\Pi_{\tau} u^\tau(t_n)\|_{l^4_\tau L^\infty} \lesssim \|\Pi_\tau u^\tau(t_{n})\|_{l^4_\tau W^{\frac d4+\delta,4}} \lesssim \|u^\tau(t_n)\|_{X^{\frac d2-\frac12+6\delta,\frac12-\delta}_\tau}\lesssim \tau^{\frac{s_0}2-3\delta-\frac d4+\frac14}
$$
for $s_0\in(\frac d2-1,\frac d2-\frac12]$. For the case $s_0>\frac d2-\frac12$, by taking $\delta<\frac15(s_0-\frac d2+\frac12)$, we can use a similar argument to get
$$
\|\Pi_\tau u^\tau(t_n)\|_{l^4_\tau L^\infty} \lesssim \|\Pi_{\tau}u^\tau(t_n)\|_{l^4_\tau W^{\frac d4+\delta,4}} \lesssim \|u^\tau(t_n)\|_{X^{s_0,\frac12-\delta}_\tau}\lesssim C_T.
$$
This yields in the case $s_0\leq\frac d2-\frac12$ that
\begin{equation}\label{ouf}
\begin{aligned}
\sum_{j=1}^5 \left\| |\Pi_{\tau }u^\tau(t_n)|^{5-j} |\Pi_\tau e_n|^j \right\|_{l^1_\tau L^2}
&\lesssim C_{T,\delta}\left(\tau^{2s_0-d+1-12\delta}
\|e_n\|_{X_\tau^{0,b_0}}+\tau^{\frac32s_0-d+1-10\delta} \|e_n\|_{X_\tau^{0,b_0}}^2 \right. \\
&\hspace{-5mm}\left. +\tau^{s_0-d+1-8\delta} \|e_n\|_{X_\tau^{0, b_0}}^3+ \tau^{ \frac{s_0}2-d+1-6\delta} \|e_n\|_{X_\tau^{0, b_0}}^4 +\tau^{-d+1-4\delta}\|e_n\|_{X_{\tau}^{0,b_0}}^5\right).
\end{aligned}
\end{equation}
The case $s_0>\frac d2-\frac12$ can be handled by similar arguments. Consequently, by combining \eqref{esten1}, \eqref{esten2}, \eqref{ouf3}, and \eqref{ouf}, we get
\begin{align*}
\|e_n\|_{X^{0, b_0}_\tau} &\leq C_T \left( ( T_1^{\varepsilon_0} + C_{T, \delta}\tau^{2s_0-d+2-13\delta-(\frac12-b_1)})\|e_n\|_{X^{0,b_0}_\tau} \right. \\
&\quad\left.{} + C_{T, s_1,\delta} (\tau^{-\frac{s_1}2} + \tau^{\frac32s_0-d+2-11\delta - (\frac12-b_1)}) \|e_n\|_{X^{0, b_0}_\tau}^2 \right. \\
&\quad\left.{} +  C_{T, s_1,\delta} ( \tau^{-s_1} +\tau^{s_0-d+2- 9\delta- (\frac12-b_1)}) \|e_n\|_{X^{0, b_0}_\tau}^3 \right. \\
&\quad\left.{} +  C_{T, \delta} \tau^{\frac{s_0}2-d+2-7\delta-(\frac12-b_1)} \|e_n\|_{X^{0, b_0}_\tau}^4+ C_{T, \delta} \tau^{-5\delta-d+2-(\frac12-b_1)}\|e_n\|_{X^{0, b_0}_\tau}^5+\tau^{\frac{s_0}2}\right).
\end{align*}
The next step is to choose the parameters appropriately. We first choose $s_1$ satisfying $\tfrac{s_0}2-\tfrac{s_1}2>0$ and $\max\big(\tfrac14,\tfrac12-\tfrac14(s_1-\frac d2+1)\big)<b^\prime<b_1$, i.e. $\frac d2+1-4b_1<\frac d2+1-4b^\prime<s_1<s_0$ (since we always have $b^\prime>\tfrac12-\tfrac14(s_0-\frac d2+1)$), and then $\delta$ sufficiently small ($s_1-\frac d2+1>100\delta$ for instance, and $\delta<\frac15(s_0-\frac d2+\frac12)$ if necessary). Since $\frac12-b_1<\tfrac14(s_1-\frac d2+1)$, we thus get
\begin{multline*}
\|e_{n}\|_{X^{0, b_0}_\tau} \leq C_T \left(( T_1^{\varepsilon_0}  + C_{T,s_1, \delta} \tau^\rho) \|e_n\|_{X^{0, b_0}_{\tau}} + C_{T, s_1,\delta}\tau^{-\frac{s_1}2} \|e_n\|_{X^{0, b_0}_\tau}^2
+ C_{T,s_1,\delta}\tau^{-s_1} \|e_n\|_{X^{0, b_0}_\tau}^3  \right.\\
\left.+ C_{T,\delta} \tau^{-\frac32s_1}\|e_n\|_{X^{0,b_0}_{\tau}}^4+ C_{T,\delta}\tau^{-2s_1}\|e_n\|_{X^{0,b_0}_\tau}^5+\tau^\frac{s_0}2\right)
\end{multline*}
for some $\rho>0$. We then choose $T_1$ and $\tau$ small enough with respect to $C_T$, so that \\$ C_{T}( T_{1}^{\varepsilon_0}  + C_{T,s_1 ,\delta}\tau^{\rho})\leq\frac12$. This yields
\begin{multline*}
\qquad\|e_n\|_{X^{0, b_0}_\tau} \leq C_{T,s_1,\delta}\tau^{-\frac{s_1}2} \|e_n\|_{X^{0, b_0}_\tau}^2 +C_{T,s_1,\delta} \tau^{-s_1} \|e_{n}\|_{X^{0, b_0}_\tau}^3 \\
+C_{T,\delta}\tau^{-\frac32s_1}\|e_n\|_{X^{0, b_0}_\tau}^4+C_{T,\delta}\tau^{-2s_1}\|e_{n}\|_{X^{0, b_0}_\tau}^5+C_T\tau^\frac{s_0}2.\qquad
\end{multline*}
Thus we finally get for $\tau$ sufficiently small that
$$
\|e_n\|_{X^{0,b_0}_\tau} \leq C_T \tau^\frac{s_0}2,
$$
which proves the desired estimate for $0 \leq n\tau \leq T_1$. Note that the choice of $T_1$ depends only on $T$, we can then reiterate the above argument on $T_1\leq n\tau\leq2T_1$ and so on, to get finally the estimate for $0\leq n\tau\leq T$.
\end{proof}

\begin{remark}\label{d6}
\rm{Similar tools can also be used for higher dimensions $d\geq6$ to prove global convergence of order one. We omit the details.}
\end{remark}

We finally estimate the global error in $L^2$. To carry out this, we first use that
\begin{align*}
\Vert u(t_n)-u_n\Vert_{L^2}&\leq\Vert u(t_n)-u^\tau(t_n)\Vert_{L^2}+\Vert u^\tau(t_n)-u_n\Vert_{L^2}\\
&\leq\Vert u-u^\tau\Vert_{L^\infty([0, T],  L^2)}+\Vert e_n\Vert_{l^\infty_{\tau}(0 \leq n\tau \leq T, L^2)}.
\end{align*}

Next, we use the embeddings \eqref{d} and \eqref{y} combined with \eqref{diffutau} and Proposition \ref{propen} to get that
$$
\Vert u(t_n)-u_n\Vert_{L^2} \leq C_T \tau^{\frac{s_0}2}.
$$
This concludes the proof of Theorem \ref{mainthm}.

\begin{figure}[hb]
\begin{center}
\subfigure[]{\includegraphics[width=0.48\textwidth]{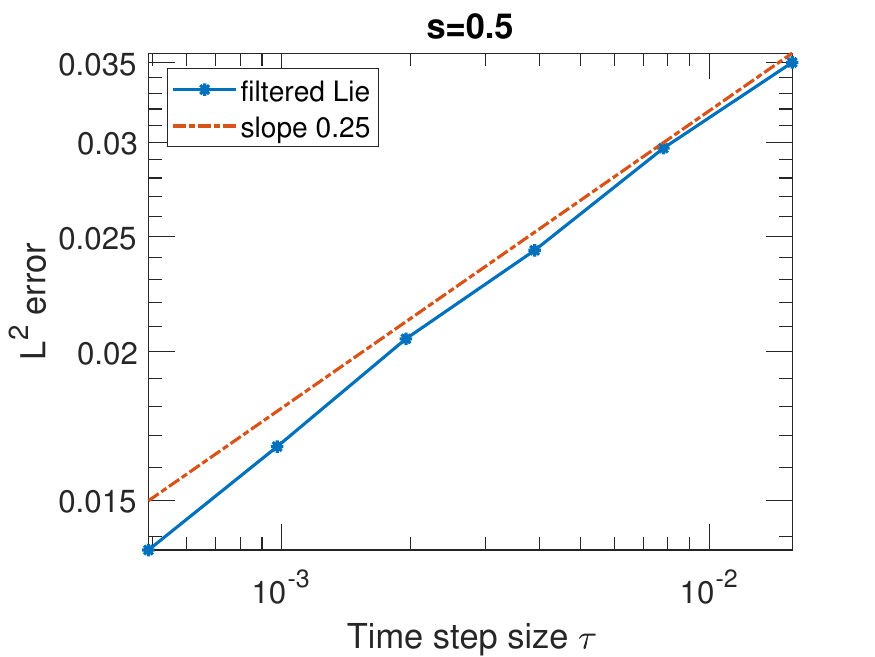}}
\subfigure[]{\includegraphics[width=0.48\textwidth]{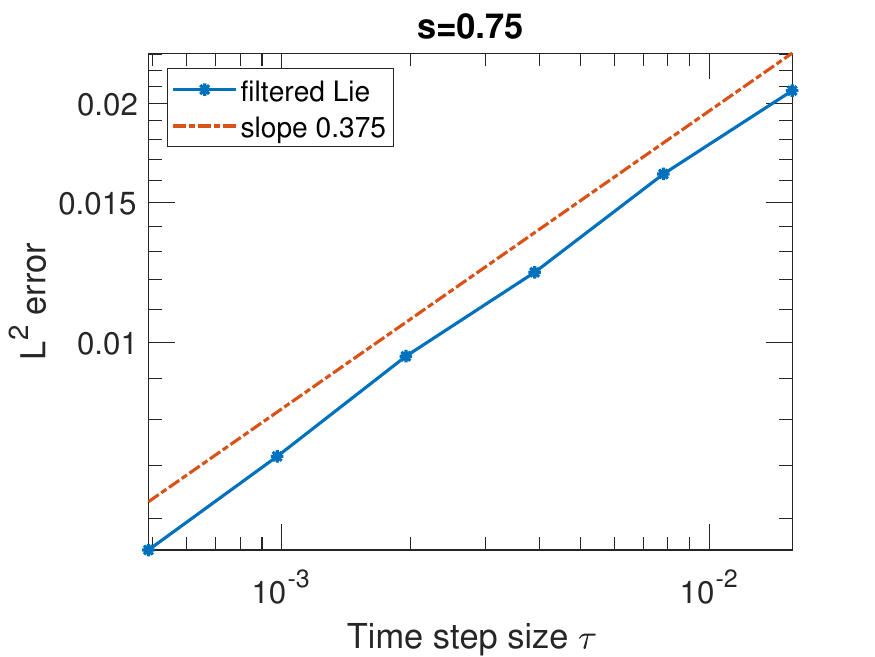}}
\subfigure[]{\includegraphics[width=0.48\textwidth]{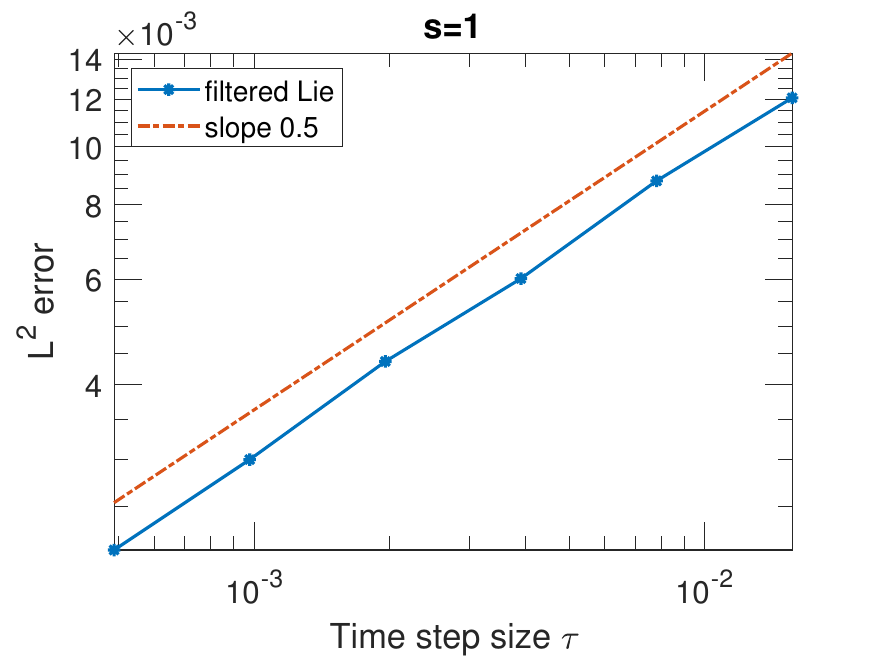}}
\subfigure[]{\includegraphics[width=0.48\textwidth]{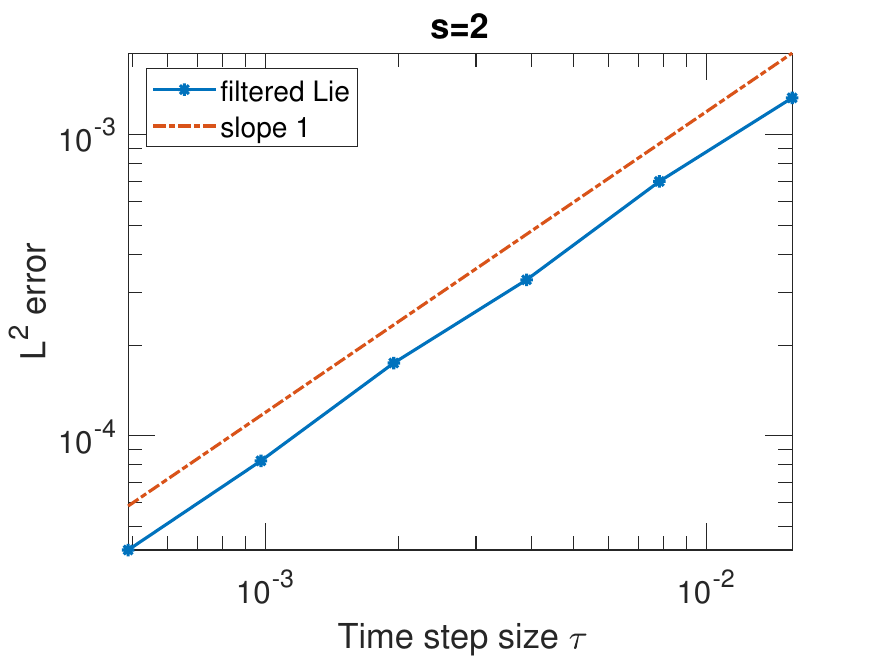}}
\end{center}
\caption{The $L^2$ error of the filtered Lie splitting scheme for the three dimensional NLS with rough initial data $u_0\in H^s(\mathbb{T}^3)$. \
(a) $s=0.5$; \quad (b) $s=0.75$; \quad (c) $s=1$; \quad (d) $s=2$.\label{figure1}}
\end{figure}

\goodbreak

\section{Numerical experiments}

To illustrate the convergence properties of the filtered Lie splitting method \eqref{lie}, we performed some numerical experiments in dimensions $d=3$ and $d=4$. As initial data, we took
$$
u_0=\sum\limits_{k\in\mathbb{Z}^d}\langle k\rangle^{-(s+\frac d2+\varepsilon)}\tilde{g}_ke^{i\langle k,x\rangle}\in H^s(\mathbb{T}^d),
$$
where $\tilde{g}_k$ are uniformly distributed random variables with values in $[-1,1]+i[-1,1]$. Here, the parameter $\varepsilon > 0$ is included to guarantee that $u_0$ lies in $H^s(\mathbb{T}^d)$. In our experiments, however, we simply chose $\varepsilon = 0$.

The case $d=3$ is shown in Figure~\ref{figure1}. For the spatial discretization of \eqref{nls}, we used a standard Fourier pseudospectral method with uniform grid size $\Delta x=0.0123$. This corresponds to the maximum Fourier mode $K=(2^9,2^9,2^9)$. In order to be able to compare the different experiments, we have normalized the initial data in the $L^2$ norm to 0.1. The reference solution was obtained with standard Lie splitting, using $2^9$ spatial points in each dimension and the time step size $\tau=2^{-16}$. In all experiments, the final time was $T=1$. Figure \ref{figure1} clearly shows the convergence rate $s/2$ for solutions in $H^s$ with $s=0.5$, $0.75$, 1, and 2. This agrees with our theoretical result given in Theorem \ref{mainthm}.

\begin{figure}[t]
\begin{center}
\subfigure[]{\includegraphics[width=0.48\textwidth]{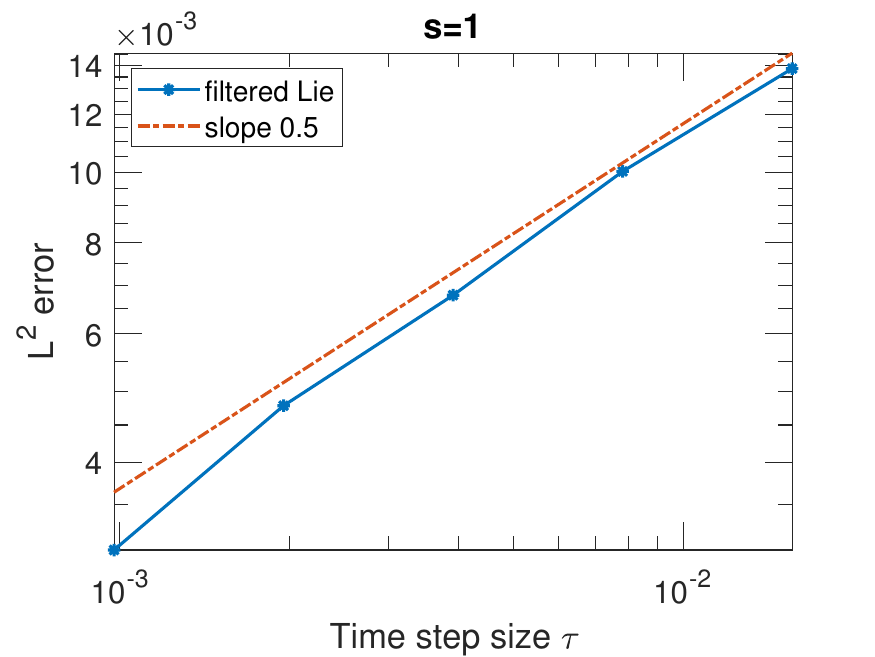}}
\subfigure[]{\includegraphics[width=0.48\textwidth]{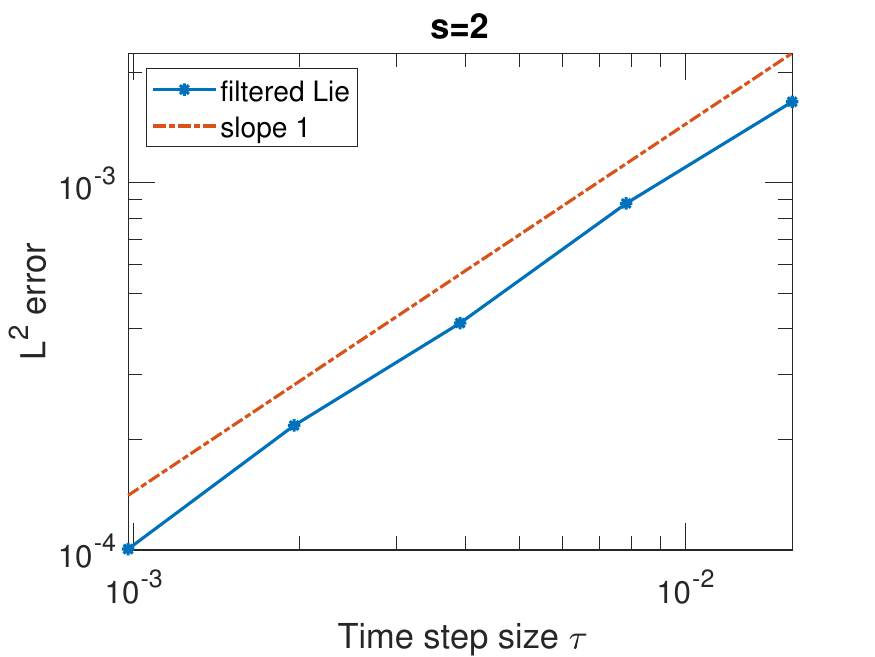}}
\end{center}
\caption{The $L^2$ error of the filtered Lie splitting scheme for the four dimensional NLS with rough initial data $u_0\in H^s(\mathbb{T}^4)$. \qquad
(a) $s=1$; \quad (b) $s=2$.\label{figure2}}.
\end{figure}

We also carried out some numerical experiments in dimension $d=4$. In this case, the grid size was chosen as $\Delta x=0.049$ which limits the Fourier modes to $K=(2^7,2^7,2^7,2^7)$. The reference solution was computed with standard Lie splitting, using $2^7$ spatial points in each dimension and the time step size $\tau=2^{-12}$. The results are shown in Figure \ref{figure2}. For $s=1$ and 2, this figure again confirms the order $s/2$ given in Theorem \ref{mainthm}.

{}
\end{document}